\documentclass[pdflatex,sn-mathphys-num]{sn-jnl}
\usepackage{amsmath,amssymb,amsthm}
\usepackage{mathrsfs}
\usepackage{tikz}
\newtheorem{theorem}{Theorem}
\newtheorem{lemma}[theorem]{Lemma}
\newtheorem{proposition}[theorem]{Proposition}
\newtheorem{corollary}[theorem]{Corollary}
\theoremstyle{remark}

\newcommand{\co}{\operatorname{co}}
\newcommand{\aff}{\operatorname{aff}}
\newcommand{\FT}{\mathcal Z}
\newcommand{\norm}[1]{\left\|#1\right\|}
\begin{document}
\title{The Fermat--Torricelli problem under isosceles orthogonality-type conditions}
\author[1]{\fnm{Zhiyao} \sur{Fang}}\email{fzy73@stu.aqnu.edu.cn}
\author*[1]{\fnm{Qi} \sur{Liu}}\email{liuq67@aqnu.edu.cn}
\affil[1]{\orgdiv{School of Mathematics and Statistics},\orgname{Anqing Normal University},\orgaddress{\city{Anqing},\postcode{246133},\state{Anhui},\country{China}}}
\abstract{The location problem studied by Weber \cite{Weber} leads to the Fermat--Torricelli problem of minimizing sums of distances. Ben\'{i}tez, Fern\'{a}ndez and Soriano \cite{BFS} proved that a real normed space of dimension at least three is an inner-product space if and only if every triple of points has a Fermat--Torricelli median in its convex hull. Pakhrou \cite{Pakhrou} showed that it suffices to test equilateral triangles. We establish a characterization using triangles $\{0,x,y\}$ satisfying
	\[
	\|x+y\|=\lambda\|x-y\|,\qquad \forall\,\lambda>2/3,
	\]
	with the parameter fixed in each characterization. The case of parameter one gives an isosceles orthogonality criterion. The proof combines radial rescaling with the transfer of planar minimality to the ambient space, without completeness or regularity assumptions on the norm.}
\maketitle

\section{Introduction}\label{sec:introduction}
For three points $a_1,a_2,a_3$ in a real normed space $X$, write
\begin{equation}\label{eq:ft-definition}
\FT_X(a_1,a_2,a_3)=\operatorname*{argmin}_{z\in X}\sum_{i=1}^3\norm{a_i-z}.
\end{equation}
Its elements are Fermat--Torricelli medians. The convex-hull property asks whether this set intersects $\co\{a_1,a_2,a_3\}$. In an inner-product space this holds for every triple. The underlying location problem goes back to Weber's formulation \cite{Weber}. Ben\'{i}tez, Fern\'{a}ndez and Soriano \cite{BFS} proved that, in real dimension at least three, the corresponding property for all triples characterizes inner-product norms. Pakhrou \cite{Pakhrou} obtained a characterization for centers defined by symmetric monotone norms and showed that testing equilateral triangles suffices. For the sum of distances, Corollary~5.6 of \cite{Pakhrou} states that a real normed space of dimension at least three is an inner-product space if and only if every equilateral triangle has an ambient Fermat--Torricelli median in its convex hull.

We restrict the test to triangles with a prescribed ratio between a median and its opposite side. For the triangle $\{0,x,y\}$, the median from the origin has length $\norm{x+y}/2$, and its opposite side has length $\norm{x-y}$. Thus prescribing their ratio amounts to prescribing the quotient $\norm{x+y}/\norm{x-y}$.

James \cite{James1945} introduced isosceles orthogonality in normed spaces by the relation $x\perp_I y$ when $\norm{x+y}=\norm{x-y}$. In a real inner-product space this is equivalent to $\langle x,y\rangle=0$. The relation is symmetric for every norm, and its homogeneity and additivity are closely tied to inner-product structure \cite{James1947,AMW}. It also describes metric bisectors: $x\perp_I y$ means that $x$ is equidistant from $y$ and $-y$. Its connections with uniqueness, convexity, geometric parameters, and orthogonality-based characterizations have been studied extensively \cite{Alonso1994,Baronti2015,Lin2005,AMW}. Related developments include Singer orthogonality and its variants, operator orthogonality and parallelism, and isosceles orthogonal triples in linear $2$-normed spaces \cite{Alonso1993,Habibzadeh2021,Bottazzi2019,Cho1992}. The case $\lambda=1$ of our criterion concerns the ambient location of Fermat--Torricelli medians for triples determined by an isosceles orthogonal pair.

\begin{theorem}\label{thm:main}
Fix $\lambda>2/3$. Let $X$ be a real normed space with $\dim X\ge3$. Then the following conditions are equivalent:
\begin{enumerate}
\item The norm of $X$ is induced by an inner product.
\item For every $x,y\in X$,
\begin{equation}\label{eq:restricted}
\norm{x+y}=\lambda\norm{x-y}
\quad\Longrightarrow\quad
\FT_X(0,x,y)\cap\co\{0,x,y\}\ne\varnothing.
\end{equation}
\end{enumerate}
\end{theorem}

The condition in \eqref{eq:restricted} reduces the family of triangles on which the convex-hull property must be tested. For $\{0,x,y\}$, it prescribes the ratio between the median from the origin and the opposite side. Theorem~\ref{thm:main} shows that testing this family at a single fixed $\lambda>2/3$ already forces the norm to come from an inner product. When $\lambda=1$, the equation is precisely isosceles orthogonality. For other values it is a prescribed-ratio condition, not an orthogonality relation in the usual sense. Indeed, in an inner-product space it is equivalent to
\[
\langle x,y\rangle=
\frac{\lambda^2-1}{2(1+\lambda^2)}
\bigl(\norm{x}^2+\norm{y}^2\bigr).
\]
For a pair not both zero, the right-hand side vanishes exactly when $\lambda=1$.

These isosceles orthogonality-type conditions suffice for the characterization because positive radial rescaling preserves Fermat--Torricelli minimality at the origin and allows the prescribed equation to be attained. To see how the restriction is removed, take a nonvertex planar median in the convex hull of an equilateral triangle and translate it to the origin. The residual vectors surround the origin, so Lemma~\ref{lem:radial-construction} provides positive radial factors that produce a test triangle. The hypothesis supplies an ambient median in its convex hull and hence in the same plane. Since the origin is still a planar median, comparison of the two minimum values makes it an ambient median as well. Inverse radial rescaling then shows that the original equilateral triangle has an ambient median in its convex hull. Pakhrou's equilateral-triangle characterization \cite[Corollary~5.6]{Pakhrou} forces the norm to be induced by an inner product. Thus the restricted family retains enough information to detect inner-product structure, including in the isosceles orthogonality case $\lambda=1$.

The lower bound on $\lambda$ ensures that the radial factors are positive. For a nonzero zero-sum triple $v_1+v_2+v_3=0$, Lemma~\ref{lem:zero-sum} gives, after relabeling, $d/\norm{v_3}\ge3/2$, where $d=\norm{v_1-v_2}$. Rescaling $v_3$ by
\[
 t=\frac12\left(\lambda\frac{d}{\norm{v_3}}-1\right)
\]
produces the required equation. The bound $\lambda>2/3$ ensures $t>0$, so the rescaled vertex remains on the same ray and Lemma~\ref{lem:radial-invariance} applies in both directions. At $\lambda=2/3$ this estimate gives only $t\ge0$; a zero factor would prevent the inverse rescaling. 

\section{Some technical lemmas}\label{sec:radial}
The subscript in $\FT_E$ specifies the space over which the sum of distances is minimized. We apply radial rescaling both in a triangle's affine plane and in the ambient space.

\begin{lemma}\label{lem:radial-invariance}
Let $E$ be a real normed space, let $u_1,u_2,u_3\in E$, and let $a_1,a_2,a_3>0$. Then
\begin{equation}\label{eq:radial-invariance}
0\in\FT_E(u_1,u_2,u_3)
\quad\Longleftrightarrow\quad
0\in\FT_E(a_1u_1,a_2u_2,a_3u_3).
\end{equation}
\end{lemma}
\begin{proof}
Assume that the origin minimizes the first sum. For any $z\in E$, choose $\delta>0$ with $\delta a_i\le1$ for every $i$. Since
\[
u_i-\delta z=(1-\delta a_i)u_i+\delta(a_iu_i-z),
\]
the triangle inequality and minimality at the origin yield
\begin{align*}
\sum_i\norm{u_i}
&\le\sum_i\norm{u_i-\delta z}\\
&\le\sum_i(1-\delta a_i)\norm{u_i}
+\delta\sum_i\norm{a_iu_i-z}.
\end{align*}
Subtracting $\sum_i\norm{u_i}$ and dividing by $\delta$ gives
\[
\sum_i\norm{a_iu_i}=\sum_i a_i\norm{u_i}
\le\sum_i\norm{a_iu_i-z}.
\]
This proves one implication. Apply it to the vectors $a_i u_i$ and the positive coefficients $a_i^{-1}$ to obtain the converse.
\end{proof}

The next estimate supplies the range of parameters in Theorem~\ref{thm:main}.
\begin{lemma}\label{lem:zero-sum}
Let $v_1,v_2,v_3$ be nonzero vectors in a normed space such that $v_1+v_2+v_3=0$. For some permutation $(i,j,k)$ of $(1,2,3)$,
\begin{equation}\label{eq:ratio-estimate}
\norm{v_i-v_j}\ge\frac32\norm{v_k}.
\end{equation}
\end{lemma}
\begin{proof}
Put $S=\sum_i\norm{v_i}$ and
$D=\norm{v_1-v_2}+\norm{v_2-v_3}+\norm{v_3-v_1}$.
For distinct $i,j,k$, the zero-sum identity gives
$3v_i=(v_i-v_j)+(v_i-v_k)$, whence
\[
3\norm{v_i}\le\norm{v_i-v_j}+\norm{v_i-v_k}.
\]
Summing gives $3S\le2D$. If \eqref{eq:ratio-estimate} failed for all three choices of the opposite index, summation would give $D<3S/2$, a contradiction.
\end{proof}

\begin{lemma}\label{lem:radial-construction}
Fix $\lambda>2/3$. Let $u_1,u_2,u_3$ be noncollinear, nonzero vectors in a real normed space, and suppose that $0\in\co\{u_1,u_2,u_3\}$. There exist $a_1,a_2,a_3>0$ such that, after relabeling the vertices, $b_i=a_iu_i$ satisfy
\begin{equation}\label{eq:constructed}
\norm{b_1+b_2-2b_3}=\lambda\norm{b_1-b_2}.
\end{equation}
\end{lemma}
\begin{proof}
First suppose the origin is in the relative interior of the triangle. Choose $\alpha_i>0$ with $\sum_i\alpha_i u_i=0$, and put $v_i=\alpha_i u_i$. By Lemma~\ref{lem:zero-sum}, relabeling gives
\[
d:=\norm{v_1-v_2}\ge\frac32\norm{v_3}.
\]
Consequently
\begin{equation}\label{eq:stretch-factor}
t=\frac12\left(\frac{\lambda d}{\norm{v_3}}-1\right)>0.
\end{equation}
Set $b_1=v_1$, $b_2=v_2$, and $b_3=tv_3$. Because $v_1+v_2=-v_3$,
\[
\norm{b_1+b_2-2b_3}
=(1+2t)\norm{v_3}
=\lambda d
=\lambda\norm{b_1-b_2}.
\]

Next suppose the origin is on the boundary of the triangle. It is not a vertex, since all $u_i$ are nonzero. Relabeling places it in the relative interior of $[u_1,u_2]$. Choose $a_1,a_2>0$ with $v=a_1u_1=-a_2u_2$, and put
\[
b_1=v,\qquad b_2=-v,\qquad
b_3=\frac{\lambda\norm{v}}{\norm{u_3}}u_3.
\]
Then $\norm{b_3}=\lambda\norm{v}$, so
\[
\norm{b_1+b_2-2b_3}=2\norm{b_3}
=2\lambda\norm{v}=\lambda\norm{b_1-b_2}.
\]

\end{proof}

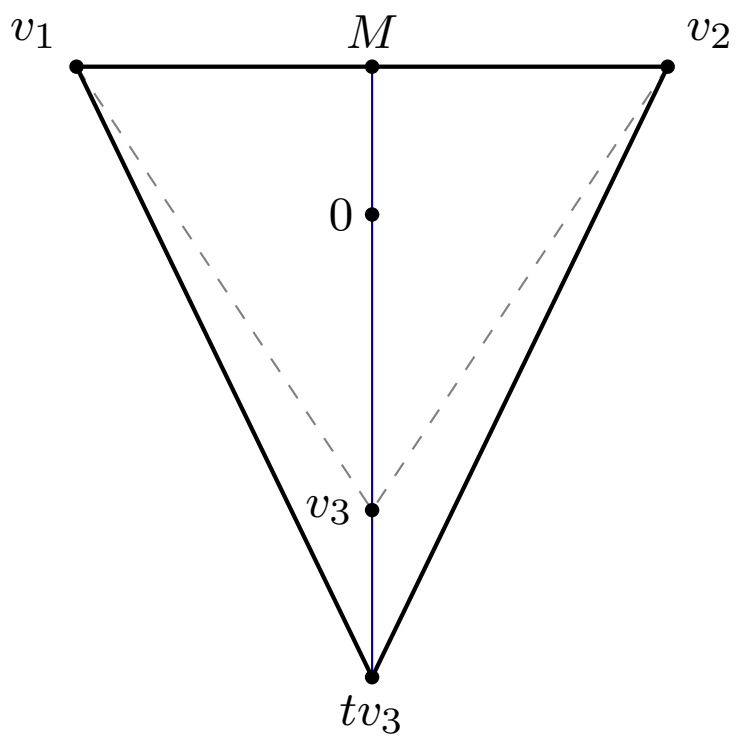
\begin{figure}[htbp]
\centering
\begin{tikzpicture}[scale=0.85, every node/.style={font=\small}]
\coordinate (A) at (-2.3,1.15);
\coordinate (B) at (2.3,1.15);
\coordinate (O) at (0,0);
\coordinate (V) at (0,-2.3);
\coordinate (C) at (0,-3.6);
\coordinate (M) at (0,1.15);
\draw[gray,dashed] (A)--(V)--(B);
\draw[thick] (A)--(B)--(C)--cycle;
\draw[blue!60!black] (M)--(C);

\foreach \p in {A,B,O,V,C,M} \fill (\p) circle (1.6pt);
\node[above left] at (A) {$v_1$};
\node[above right] at (B) {$v_2$};
\node[above] at (M) {$M$};
\node[left] at (O) {$0$};
\node[left] at (V) {$v_3$};
\node[below] at (C) {$tv_3$};
\end{tikzpicture}
\caption{Radial rescaling $v_3\mapsto tv_3$. Dashed edges show the original triangle; solid edges show the rescaled triangle. The median joins $tv_3$ to $M=(v_1+v_2)/2$.}\label{fig:radial}
\end{figure}

The interior construction in Lemma~\ref{lem:radial-construction} is shown in Figure~\ref{fig:radial}. The zero-sum relation places the midpoint $M=(v_1+v_2)/2=-v_3/2$ on the line through $0$ and $v_3$. Moving the third vertex to $tv_3$ therefore prescribes the ratio of the median to the fixed opposite side.

\section{Proof of the main Theorem}\label{sec:transfer}
The following proposition transfers a planar minimizer to the ambient space.
\begin{proposition}\label{prop:transfer}
Suppose that $X$ satisfies \eqref{eq:restricted} for a fixed $\lambda>2/3$. Let $A_1,A_2,A_3$ be noncollinear points and let $P=\aff\{A_1,A_2,A_3\}$. If
\[
m\in\FT_P(A_1,A_2,A_3)\cap\co\{A_1,A_2,A_3\},
\qquad m\notin\{A_1,A_2,A_3\},
\]
then $m\in\FT_X(A_1,A_2,A_3)$.
\end{proposition}
\begin{proof}
Translate $m$ to the origin and write $Y=P-m$ and $u_i=A_i-m$. The assumptions give nonzero $u_i$ with
\[
0\in\FT_Y(u_1,u_2,u_3),\qquad 0\in\co\{u_1,u_2,u_3\}.
\]
Choose $a_i>0$ by Lemma~\ref{lem:radial-construction}, and put $b_i=a_i u_i$, so that \eqref{eq:constructed} holds. Lemma~\ref{lem:radial-invariance} in $Y$ gives $0\in\FT_Y(b_1,b_2,b_3)$.

Apply \eqref{eq:restricted} to $x=b_1-b_3$ and $y=b_2-b_3$. Translation gives a point
\[
q\in\FT_X(b_1,b_2,b_3)\cap\co\{b_1,b_2,b_3\}.
\]
In particular $q\in Y$. Planar minimality at $0$ and ambient minimality at $q$ imply
\[
\sum_i\norm{b_i}\le\sum_i\norm{b_i-q}\le\sum_i\norm{b_i}.
\]
Thus $0\in\FT_X(b_1,b_2,b_3)$. Apply Lemma~\ref{lem:radial-invariance} in $X$, with the inverse radial factors, to obtain $0\in\FT_X(u_1,u_2,u_3)$. Translating back proves the assertion.
\end{proof}

We recall the planar result of Wendell and Hurter as stated in \cite[Theorem~3.23]{MSW}: every finite subset of a normed plane has a Fermat--Torricelli median in its convex hull. Planar minimizers exist because the distance sum is continuous and coercive in finite dimension. For equilateral triples, a minimizer can be chosen away from the vertices.

\begin{lemma}\label{lem:equilateral-planar}
Let $A_1,A_2,A_3$ form an equilateral triangle in a normed plane $P$. There is a point
\[
m\in\FT_P(A_1,A_2,A_3)\cap\co\{A_1,A_2,A_3\}
\]
which is different from all three vertices.
\end{lemma}
\begin{proof}
Let $d>0$ be the common side length, put $F(z)=\sum_i\norm{A_i-z}$, and write $\mu=\min_{z\in P}F(z)$. The planar convex-hull result supplies a minimizer in the triangle. Each vertex has value $2d$, so $\mu\le2d$. If $\mu<2d$, that minimizer is not a vertex. If $\mu=2d$, all three vertices minimize $F$. Its convexity gives
\[
F\left(\frac{A_1+A_2+A_3}{3}\right)
\le\frac13\sum_{i=1}^3 F(A_i)=2d=\mu.
\]
Hence the centroid is a minimizer. A triple of distinct equidistant points cannot be collinear: on a line, the distance between the two extreme points is the sum of the other two distances. Its centroid therefore differs from each vertex.
\end{proof}

\begin{proof}[Proof of Theorem~\ref{thm:main}]
Assume condition~2. For any equilateral triangle in $X$, Lemma~\ref{lem:equilateral-planar} gives a nonvertex minimizer in its convex hull when minimization is restricted to its affine plane. Proposition~\ref{prop:transfer} makes this point an ambient minimizer. Thus every equilateral triangle has an ambient Fermat--Torricelli median in its convex hull. By \cite[Corollary~5.6]{Pakhrou}, $X$ is an inner-product space.

Conversely, suppose the norm comes from an inner product. For any triple, let $C=\co\{A_1,A_2,A_3\}$. The continuous distance sum attains its minimum at some $m\in C$, since $C$ is compact. For $z\notin C$, its distance to $C$ attains a minimum at $c\in C$. Differentiating the squared distance along $c+s(a-c)$ for $a\in C$ and $s\downarrow0$ gives
\[
\langle z-c,a-c\rangle\le0.
\]
It follows for each $i$ that
\[
\norm{A_i-z}^2
=\norm{A_i-c}^2+\norm{z-c}^2
-2\langle A_i-c,z-c\rangle
\ge\norm{A_i-c}^2.
\]
Thus the distance sum at $z$ is at least its value at $c$, which in turn is at least its value at $m$. The same conclusion for $z\in C$ follows from the definition of $m$. Therefore $m$ is an ambient median. The convex-hull property follows for every triple and hence for those in \eqref{eq:restricted}. The argument uses compactness of $C$ and applies equally to incomplete inner-product spaces.
\end{proof}

\label{sec:consequences}
\begin{corollary}\label{cor:isosceles}
Let $X$ be a real normed space with $\dim X\ge3$. Its norm is induced by an inner product if and only if
\[
x\perp_I y\quad\Longrightarrow\quad
\FT_X(0,x,y)\cap\co\{0,x,y\}\ne\varnothing
\qquad(x,y\in X).
\]
\end{corollary}
\begin{proof}
Take $\lambda=1$ in Theorem~\ref{thm:main}.
\end{proof}

For $\lambda=1$, the median from the origin has half the length of the opposite side. This is distinct from the equal-side condition $\norm{x}=\norm{y}$ for the triangle $\{0,x,y\}$.

\begin{corollary}\label{cor:median-side}
Fix $r>1/3$. A real normed space $X$ of dimension at least three is an inner-product space if and only if every nondegenerate triangle having a median of length $r$ times its opposite side has an ambient Fermat--Torricelli median in its convex hull.
\end{corollary}
\begin{proof}
At the vertex $A_3$, put $x=A_1-A_3$ and $y=A_2-A_3$. The prescribed geometric ratio is
\[
\frac{\norm{(A_1+A_2)/2-A_3}}{\norm{A_1-A_2}}
=\frac{\norm{x+y}}{2\norm{x-y}}=r.
\]
Apply Theorem~\ref{thm:main} with $\lambda=2r$. Its proof of the nontrivial implication only uses nondegenerate test triangles: in the interior case, the relation $v_1+v_2+v_3=0$ has positive coefficients, and replacing $v_3$ by $tv_3$ with $t>0$ preserves noncollinearity. Indeed, an affine dependence of $v_1,v_2,tv_3$ would require coefficients proportional to $(t,t,1)$, whose sum is nonzero. In the boundary case, $b_3$ lies off the line through $v$ and $-v$. Thus the nondegenerate restriction is sufficient.
\end{proof}

The isosceles criterion fails in $\ell_1^3$. Indeed, take $X=(\mathbb R^3,\norm{\cdot}_1)$ and
\[
x=(1,1,1),\qquad y=(1,-1,0).
\]
Then $\norm{x+y}_1=\norm{x-y}_1=3$. The unique ambient minimizer for $\{0,x,y\}$ is $m=(1,0,0)$: minimization separates into three scalar sums of absolute values, whose unique minimizers are the respective coordinate medians. If $m=ax+by$ with $a,b\ge0$ and $a+b\le1$, the third coordinate gives $a=0$, and the second then gives $b=0$, contradicting the first. Hence $m$ is outside the triangle. 

The dimensional assumption is necessary. In dimension two the planar convex-hull theorem gives the required property for every triple and every parameter, including for non-Euclidean norms. The bound $\lambda>2/3$ ensures that the radial factor in \eqref{eq:stretch-factor} is positive by \eqref{eq:ratio-estimate}. 

\section*{Acknowledgements}

Thanks to all the members of the Functional Analysis Research team of the College of Mathematics and Statistics of Anqing Normal University for their discussion and correction of the difficulties and errors encountered in this paper. The authors cordially thank Professor Hai Zhang for his encouraging advice and the first author would like to acknowledge his strong support for undergraduate students in carrying out basic mathematics research.

\section*{Declarations}

\subsection*{Competing Interests}
The authors declare that there are no competing interests.

\subsection*{Funding Information}
Not Applicable.

\subsection*{Author contribution}
These authors contributed equally to this work. All authors have read and agreed to the published version of the manuscript.

\subsection*{Data Availability Statement}
Not Applicable.

\end{document}